\documentclass[12pt,a4paper,reqno]{amsart}
\usepackage{fullpage}
\usepackage{verbatim}
\usepackage{hyperref}
\usepackage{pslatex} 

\usepackage[full]{textcomp}
\usepackage[osf]{newtxtext} 
\usepackage{cabin} 
\usepackage[bigdelims,vvarbb]{newtxmath}  
\usepackage[cal=boondoxo]{mathalfa}   
 
\usepackage{microtype}  

\newcommand{\eps}{\varepsilon}
\newcommand{\sameorder}{\asymp}
\newcommand{\dx}{\mathrm{d}}
\newcommand{\I}{\mathcal{I}}
\newcommand{\J}{\mathcal{J}}
\newcommand{\M}{\mathcal{M}}
\newcommand{\C}{\mathbb{C}}

\newcommand{\R}{\mathbb{R}}
\newcommand{\Y}{\mathbb{Y}}

\newcommand{\Stilde}{\widetilde{S}}
\newcommand{\Odip}[2]{\mathcal{O}_{#1}\left(#2\right)}
\newcommand{\Odipg}[2]{\mathcal{O}_{#1}\bigl(#2\bigr)}
\newcommand{\Odipm}[2]{\mathcal{O}_{#1} (#2)}
\newcommand{\Odi}[1]{\Odip{}{#1}}
\newcommand{\Odim}[1]{\Odipm{}{#1}}

\newtheorem{Theorem}{Theorem}[section]
\newtheorem{Lemma}{Lemma}[section]

\allowdisplaybreaks

\title[{A Ces\`aro Average of generalised  Hardy-Littlewood numbers}]{A Ces\`aro Average of generalised\\ Hardy-Littlewood numbers}

\author{Alessandro Languasco \lowercase{and} Alessandro Zaccagnini}


\begin{document}

\begin{abstract}
We continue our recent work on additive problems with prime summands:
we already studied the \emph{average} number of representations of an
integer as a sum of two primes, and also considered individual integers.
Furthermore, we dealt with representations of integers as sums of powers
of prime numbers.
In this paper, we study a Ces\`aro weighted partial \emph{explicit} formula for
generalised Hardy-Littlewood numbers (integers that can be written as
a sum of a prime power and a square) thus extending and improving our
earlier results.
\end{abstract}

\maketitle

\section{Introduction}

The problem of counting the number of representations of an integer as
a sum of some fixed powers of primes, and its variants where some
primes are replaced by powers of integers, has received much attention
in the last decades.
The goal, that has been attained only in part, is to obtain an
asymptotic formula for the number of such representations, which is
valid for large integers that satisfy some necessary congruence
conditions, as in the binary and ternary Goldbach problems.
In some cases, conditional results are obtained, that is, it is
necessary to assume the truth of some hitherto unproved conjecture like
the Riemann Hypothesis.
In our previous paper \cite{LanguascoZ2012b} we considered the problem
of representing a large integer of suitable parity as a sum of
$j \ge 5$ primes, assuming the truth of the Generalised Riemann
Hypothesis, and we obtained an individual asymptotic formula with a
main term of the expected order of magnitude, and a lower order term
which depends explicitly on non-trivial zeros of relevant Dirichlet
$L$-functions.
The corresponding problem with a smaller number of summands is harder,
and it is convenient to study the \emph{average} number of
representations; in fact, assuming the Riemann Hypothesis, in our paper
\cite{LanguascoZ2012a} we gave one such result for the standard
Goldbach problem where $j = 2$.
In this case, averaging has the effect of making the zeros of
$L$-functions irrelevant, except for the Riemann $\zeta$-function
itself, and in fact the development contains a main term and a smaller
term which depends on the non-trivial zeros of the $\zeta$-function.

In the paper \cite{LanguascoZ2015a} we introduced a Ces\`aro-Riesz
weight in the summation: the presence of a smooth weight in place,
essentially, of the characteristic function of the interval where we
are averaging, leads to the possibility of giving a development into
several terms, of decreasing order of magnitude, depending on the
zeros of the Riemann $\zeta$-function, sometimes in pairs.
This weight also enabled us to remove the necessity of assuming
the Riemann Hypothesis.
The results of this paper have been generalised and improved in
\cite{LanguascoZ2017c}, where we treated the average number of
representations of an integer in the form $p_1^{\ell_1} + p_2^{\ell_2}$,
where $\ell_1$ and $\ell_2$ are fixed positive integers and $p_1$
and $p_2$ are prime numbers.
The results in \cite{LanguascoZ2015a} have been recently extended by
Goldston and Yang \cite{GoldstonY2017} and by Br\"udern, Kaczorowski
and Perelli \cite{BrudernKP2017}.

As remarked above, we also considered a mixed binary problem with a
prime and the square of an integer, the so-called Hardy-Littlewood
numbers, in \cite{LanguascoZ2013a}.
Similar problems have been studied by Cantarini in \cite{Cantarini2017}
and \cite{Cantarini2018}.
Our task here is to extend and improve our earlier results on weighted
averages.
We let $\ell\ge 1$ be an integer and set
\begin{equation}
\label{r2-def}
   r_{\ell,2}(n)  =  \sum_{m_1^\ell + m_2^2 = n} \Lambda(m_1),
\end{equation}
where $\Lambda$ is the usual von Mangoldt-function.
Our main goal is to give a multi-term development for
\begin{equation}
\label{Rk-def}
  R_k(N)
  =
  R_k(N; \ell)
  =
  \sum_{n \le N} r_{\ell,2}(n) \frac{(1 - n/N)^k}{\Gamma(k + 1)},
\end{equation}
where $k > 0$.
We introduce the following abbreviations for the terms of the
development:
\begin{align*}
  \M_{1, \ell, k}(N)
  &=
  \pi^{1/2} \frac{\Gamma(1/\ell)}{2\ell}
  \frac{N^{1/2+1/\ell}}{\Gamma(k + 3/2+1/\ell)}
  -
  \frac{\Gamma(1/\ell)}{2\ell}
  \frac{N^{1/\ell}}{\Gamma(k + 1 + 1/\ell)}, \\
  \M_{2, \ell, k}(N)
  &=
  -
  \frac{\pi^{1/2}}{2\ell}
  \sum_{\rho}
  \frac{\Gamma(\rho/\ell)}
   {\Gamma(k + 3 / 2 + \rho/\ell)}
   N^{\rho/\ell + 1/2}, \\
  \M_{3, \ell, k}(N)
  &=
  \frac{1}{2\ell}
  \sum_{\rho}
    \frac{\Gamma(\rho/\ell)}{\Gamma(k + 1 + \rho/\ell)} N^{\rho/\ell}, \\
  \M_{4, \ell, k}(N)
  &=
  -
  \frac {\pi^{1/2}\log(2\pi)}{2\Gamma(k+3/2)} N^{1/2}, \\
  \M_{5, \ell, k}(N)
  &=
  \frac{N^{1 / 4 -k / 2 + 1/(2\ell)}}{\pi^{k + 1/\ell}}
  \frac{\Gamma(1/\ell)}{\ell}
  \sum_{j \ge 1}
    \frac{J_{k + 1 / 2 + 1/ \ell} (2 \pi j N^{1 / 2})}
         {j^{k + 1 / 2 + 1/ \ell}}, \\
  \M_{6, \ell, k}(N)
  &=
  -
  \frac{N^{1 / 4 -k / 2}}{\pi^k}
  \sum_{\rho} \Gamma\Bigl(\frac{\rho}{\ell}\Bigr)
    \frac{N^{\rho / (2\ell)}}{\pi^{\rho/\ell}}
    \sum_{j \ge 1}
      \frac{J_{k + 1 / 2 + \rho/\ell} (2 \pi j N^{1 / 2})}
           {j^{k + 1 / 2 + \rho/\ell}}, \\
  \M_{7, \ell, k}(N)
  &=
- \frac{\log(2\pi)}{\pi^{k}} N^{1 / 4 -k / 2}
     \sum_{j \ge 1}
      \frac{J_{k + 1 / 2} (2 \pi j N^{1 / 2})}
           {j^{k + 1 / 2}}.
\end{align*}
Here $\rho$ runs over the non-trivial zeros of the Riemann
zeta-function $\zeta(s)$, $\Gamma$ is Euler's function and
$J_{\nu} (u)$ denotes the Bessel function of complex order $\nu$ and
real argument $u$.
The main result of the paper is the following theorem.

\begin{Theorem}
\label{Cesaro-gen-HL-average}
Let $\ell\ge 1$ be an integer and $N$ be a sufficiently large integer.
For $k > 1$ we have
\[
  \sum_{n \le N} r_{\ell,2}(n) \frac{(1 - n/N)^k}{\Gamma(k + 1)}
  =
  \sum_{j=1}^7 \M_{j,\ell,k}(N)
  +
  \Odip{k,\ell}{1}.
\]
\end{Theorem}

Clearly, depending on the size of $\ell$, some of the previously listed
terms can be included in the error term.
Theorem \ref{Cesaro-gen-HL-average} generalises and improves our
Theorem 1 in \cite{LanguascoZ2013a}, which corresponds to the case
$\ell=1$, where the error term there should be read as
$\Odipm{k}{N^{1/2}}$; see Theorem 2.3 of \cite{Languasco2016a}.
In fact, in this case we are now able to detect the terms $\M_{4,1,k}$
and $\M_{7,1,k}$.

The basic strategy of the proof depends on the modern version of a
classical formula due to Laplace \cite{Laplace1812}, namely
\begin{equation}
\label{Laplace-transf}
  \frac 1{2 \pi i}
  \int_{(a)} v^{-s} e^v \, \dx v
  =
  \frac1{\Gamma(s)},
\end{equation}
where $\Re(s) > 0$ and $a > 0$; see Formula 5.4(1) on page 238
of Erd\'elyi et al, \cite{ErdelyiMOT1954a}.
Using a suitable form of this transform, which we describe in
\S\ref{settings}, we obtain the fundamental relation for the method,
viz.
\begin{equation}
\label{main-form-omega}
  N^k
  R_k(N)
  =
  \sum_{n \le N}  r_{\ell,2}(n) \frac{(N - n)^k}{\Gamma(k + 1)}
  =
  \frac 1{2 \pi i}
  \int_{(a)} e^{N z} z^{- k - 1} \Stilde_\ell(z) \omega_2(z) \, \dx z,
\end{equation}
where
\begin{equation}
\label{Stilde-omega-def}
  \Stilde_\ell(z)
  =
  \sum_{m \ge 1} \Lambda(m) e^{- m^\ell z}
  \quad\text{and}\quad
  \omega_2(z)
  =
  \sum_{m \ge 1} e^{-m^2 z}
\end{equation}
are the exponential sums that embody the properties of the $\ell$-th
powers of primes, and of the perfect squares, respectively.
Here we need $k > 0$, and consider the complex variable $z = a + i y$
with $a > 0$.

The basic facts that we need are the ``explicit formula'' for
$\Stilde_\ell(z)$, that is, its development as a main term and a
secondary term which is a sum over non-trivial zeros of the Riemann
$\zeta$-function with a very small error, as in \eqref{def-Sl} below,
and the simple connection of $\omega_2(z)$ with
$\theta(z) = \sum_{m = -\infty}^{+\infty} e^{-m^2 z}$, since
$\theta(z) = 1 + 2 \omega_2(z)$.
Now, we recall that $\theta$ satisfies the functional equation
\eqref{func-eq-theta}.
We plug these relations into the right-hand side of
\eqref{main-form-omega}, and exchange summation over zeros with
vertical integration, obtaining formally the development in
Theorem~\ref{Cesaro-gen-HL-average}.
The Bessel functions in $\M_{5,\ell,k}$, $\M_{6,\ell,k}$ and
$\M_{7,\ell,k}$ arise from the ``modular'' terms in the functional
equation of $\theta$.

Of course, we need to prove that the exchange referred to above is
legitimate, and that the error term arising from the approximation of
the exponential sum $\Stilde_\ell(z)$ in \eqref{def-Sl} is small.

Summing up, as in \cite{LanguascoZ2017c} we combine the approach with line
integrals with the classical methods dealing with infinite sums over
primes, exploited by Hardy and Littlewood (see \cite{HardyL1916} and
\cite{HardyL1923}) and by Linnik \cite{Linnik1946}.
The main difficulty here is, as in \cite{LanguascoZ2013a}, that the
problem naturally involves the modular relation for the complex theta
function~\eqref{func-eq-theta}.
The presence of the Bessel functions
in our statement strictly depends on such modularity relation.
It is worth mentioning that it is not clear how to get such ``modular'' terms
using the finite sums approach for the function
$ r_{\ell,2}(n)$.
The previously mentioned improvement we get in Theorem
\ref{Cesaro-gen-HL-average} follows using Lemma \ref{Linnik-lemma-gen}
below, which is proved in \cite{LanguascoZ2017c}.

\section{Settings}
\label{settings}

As we mentioned in the previous section, we will need the general case
of \eqref{Laplace-transf}, which can be found in de Azevedo Pribitkin
\cite{Azevedo2002}, formulae (8) and (9).
More precisely, we have
\begin{equation}
\label{Laplace-eq-1}
  \frac1{2 \pi}
  \int_{\R} \frac{e^{i D u}}{(a + i u)^s} \, \dx u
  =
  \begin{cases}
    \dfrac{D^{s - 1} e^{- a D}}{\Gamma(s)}
    & \text{if $D > 0$,} \\
    0
    & \text{if $D < 0$,}
  \end{cases}
\end{equation}
which is valid for $\sigma = \Re(s) > 0$ and $a \in \C$ with
$\Re(a) > 0$, and
\begin{equation}
\label{Laplace-eq-2}
  \frac1{2 \pi}
  \int_{\R} \frac 1{(a + i u)^s} \, \dx u
  =
  \begin{cases}
    0     & \text{if $\Re(s) > 1$,} \\
    1 / 2 & \text{if $s = 1$,}
  \end{cases}
\end{equation}
for $a \in \C$ with $\Re(a) > 0$.
Formulae \eqref{Laplace-eq-1}-\eqref{Laplace-eq-2} actually enable us
to write averages of arithmetical functions by means of line integrals
as we will see below.

We will also need  Bessel functions of complex order $\nu$ 
and real argument $u$.
For their definition and main properties we
refer to Watson \cite{Watson1966}.
In particular, equation (8) on page 177 gives the
Sonine representation:
\begin{equation}
  \label{Bessel-def}
  J_\nu(u)
  :=
  \frac{(u / 2)^\nu}{2 \pi i}
  \int_{(a)} s^{- \nu - 1} e^s e^{- u^2 / 4 s} \, \dx s,
\end{equation}
where $a > 0$ and $u,\nu \in \C$ with $\Re(\nu) > -1$.
We will also use the Poisson integral formula 
\begin{equation}
\label{Poisson-int-rep}
J_\nu(u)
: =
\frac{2(u/2)^{\nu}}{\pi^{1/2}\Gamma(\nu+1/2)}
\int_{0}^{1} (1-t^2)^{\nu-1/2} \cos (ut)\ \dx t,
\end{equation}
which holds for $\Re (\nu) > -1/2$ and $u\in \C$. 
(See eq.~(3) on page 48 of \cite{Watson1966}.)
An asymptotic estimate we will need is  
\begin{equation}
\label{Lebedev-asymp}
J_\nu(u)
=
\Bigl(\frac{2}{\pi u}\Bigr)^{1/2}
\cos \Bigl(u -\frac{\pi \nu}{2} -\frac{\pi}{4}\Bigr)
+
\Odipg{\vert \nu \vert}{u^{-5/2}},
\end{equation}
which follows from eq.~(1) on page 199 of Watson \cite{Watson1966}.

From now on we assume that $k > 0$.
We recall the definitions~\eqref{Stilde-omega-def}, where
$z = a + i y$ with $a > 0$.
We also recall that
$\theta(z) = \sum_{m = -\infty}^{+\infty} e^{-m^2 z}$ satisfies the
functional equation
\begin{equation}
\label{func-eq-theta}
  \theta(z)
  =
  \Bigl( \frac \pi z \Bigr)^{1/2}
  \theta\Bigl( \frac{\pi^2} z \Bigr)\quad
  \text{for} \ \Re (z)>0.
\end{equation}
See, e.g., Proposition VI.4.3 of Freitag and Busam
\cite[page 340]{FreitagB2009}.
Since $\theta(z) = 1 + 2 \omega_2(z)$, we immediately get 
\begin{equation}
\label{func-eq-omega}
  \omega_2(z)
  =
  \frac 12
  \Bigl( \frac \pi z \Bigr)^{1 / 2}
  -
  \frac12
  +
  \Bigl( \frac \pi z \Bigr)^{1 / 2}
  \omega_2 \Bigl( \frac {\pi^2} z \Bigr)
  \quad
    \textrm{for} \ \Re (z)>0.
\end{equation}
Recalling \eqref{r2-def}, we can write
\[
  \Stilde_\ell(z) \omega_2(z)
  =
  \sum_{m_1 \ge 1} \sum_{m_2 \ge 1} \Lambda(m_1) e^{-(m_1^\ell + m_2^2) z}
  =
  \sum_{n \ge 1}  r_{\ell,2}(n) e^{- n z}
\]
and, by \eqref{Laplace-eq-1}-\eqref{Laplace-eq-2}, we see  that
\begin{equation}
\label{first-step}
\sum_{n \le N}  r_{\ell,2}(n) \frac{(N - n)^k}{\Gamma(k + 1)}
 =
\sum_{n \ge 1}  r_{\ell,2}(n) 
\Bigl(
 \frac{1}{2 \pi i}
  \int_{(a)} e^{(N- n)z} z^{- k - 1}  \, \dx z
  \Bigr).
\end{equation}

Our first goal is to exchange the series with the line integral in
\eqref{first-step}. To do so
we have to recall that the Prime Number Theorem (PNT) implies the statement
\begin{equation*}
  \Stilde_\ell(a)
  \sim
 \frac{\Gamma(1/\ell)}{\ell a^{1/\ell}}
  \qquad\text{for $a \to 0+$.}
\end{equation*}

In fact, by a straightforward application of the partial summation
formula we see that
\[
  \Stilde_\ell(a)
  =
  \ell a \int_0^{+\infty} \psi(t) t^{\ell - 1} e^{-a t^\ell} \, \dx t,
\]
where $\psi$ is the standard Chebyshev function.
It is now convenient to split the integration range at
$t_0 = a^{-1 / (2 \ell)}$.
We use the weak upper bound $\psi(t) \ll t$ on $[0, t_0]$,
recalling that $\psi(t) = 0$ for $t < 2$.
Since $t^{\ell} e^{-a t^\ell} \le (a e)^{-1}$ for all $t \ge 0$, we
immediately see that the contribution of this range to $\Stilde_\ell$
is $\ll \ell a^{-1 / (2 \ell)}$.
According to a weak form of the PNT, we have
$\psi(t) = t + \Odip{A}{t (\log t)^{-A}}$ for any fixed $A > 0$.
Hence, completing the missing range, performing the obvious change
of variables and using the same bounds as above when needed, we have
\begin{align*}
  \int_{t_0}^{+\infty} \psi(t) t^{\ell - 1} e^{-a t^\ell} \, \dx t
  &=
  \int_0^{+\infty} t^{\ell} e^{-a t^\ell} \, \dx t
  +
  \Odi{ \int_0^{t_0} t^{\ell} e^{-a t^\ell} \, \dx t }
  +
  \Odi{
    \int_{t_0}^{+\infty} \frac{t^{\ell} e^{-a t^\ell}}{(\log t)^A} \, \dx t} \\
  &=
  \frac{\Gamma(1 + 1 / \ell)}{\ell a^{1 + 1 / \ell}}
  +
  \Odim{t_0 a^{-1}}
  +
  \Odi{(\log t_0)^{-A} \int_0^{+\infty} t^{\ell} e^{-a t^\ell} \, \dx t } \\
  &=
  \frac{\Gamma(1 + 1 / \ell)}{\ell a^{1 + 1 / \ell}}
  +
  \Odip{\ell}{a^{- 1 - 1 / \ell} \bigl( \log (1/a) \bigr)^{-A}}.
\end{align*}
The final result then follows by recalling that $\Gamma(s+1)=s\Gamma(s)$.
 
We will also use the inequality
\begin{equation}
\label{omega-estim}
  \vert \omega_2(z)\vert 
  \le
  \omega_2(a)  
  \le
  \int_{0}^{\infty} e^{-at^{2}} \dx t
   \le
   a^{- 1 / 2}
  \int_{0}^{\infty} e^{-v^{2}} \dx v
  \ll
  a^{- 1 / 2},
\end{equation}
from which we immediately get
\begin{align*}
  \sum_{n \ge 1} \bigl\vert   r_{\ell,2}(n) e^{- n z} \bigr\vert 
  &=
  \sum_{n \ge 2}    r_{\ell,2}(n) e^{- n a}  
=
  \Stilde_\ell(a) \omega_2(a)
  \ll_\ell
  a^{-1/\ell - 1 / 2}.
\end{align*}
Taking into account the estimates
\begin{equation}
\label{z^-1}
  \vert z \vert^{-1}
  \sameorder
  \begin{cases}
    a^{-1}   &\text{if $\vert y \vert \le a$,} \\
    \vert y \vert^{-1} &\text{if $\vert y \vert \ge a$,}
  \end{cases}
\end{equation}
 where $f\sameorder g$ means $g \ll f \ll g$, and
\[
  \vert e^{N z} z^{- k - 1}\vert 
  \sameorder
   e^{N a}
  \begin{cases}
    a^{- k - 1} &\text{if $\vert y \vert \le a$,} \\
    \vert y \vert^{- k - 1} &\text{if $\vert y \vert \ge a$,}
  \end{cases}
\]
we have
\begin{align*}
  \int_{(a)} \vert e^{N z} z^{- k - 1}\vert  \,
    \vert 
     \Stilde_\ell(z) \omega_2(z)
    \vert  \, \vert \dx z \vert
  &\ll_\ell
  a^{-1/\ell - 1 / 2} e^{N a}
  \Bigl(
    \int_{-a}^a a^{- k - 1} \, \dx y
    +
    2
    \int_a^{+\infty} y^{- k - 1} \, \dx y
  \Bigr) \\
  &\ll_\ell
   a^{-1/\ell - 1 / 2} e^{N a}
  \Bigl( a^{-k} + \frac{a^{-k}}k \Bigr).
\end{align*}
The last estimate is valid only if $k > 0$.
So, for $k > 0$, we can exchange the line integral with the sum over $n$
in \eqref{first-step}, thus proving~\eqref{main-form-omega}.

\section{Inserting zeros and modularity}

We need $k>  1/2$ in this section.
The treatment of the integral on the right-hand side
of \eqref{main-form-omega} requires Lemma \ref{Linnik-lemma-gen}.
We split $\Stilde_{\ell}(z)$ according to its statement as
$\mathcal{S}_{\ell}(z) + E(a,z,\ell)$ where $E$ satisfies the
bound in \eqref{expl-form-err-term-strong} and
\begin{equation}
\label{def-Sl}  
\mathcal{S}_{\ell}(z)
:= 
\frac{\Gamma(1/\ell)}{\ell z^{1/\ell}}
- 
\frac{1}{\ell}\sum_{\rho}z^{-\rho/\ell}\Gamma\Bigl(\frac{\rho}{\ell}\Bigr) 
- \log (2\pi),
\end{equation}
where $\rho = \beta + i\gamma$ runs over the non-trivial zeros of
$\zeta(s)$.
Formula \eqref{main-form-omega} becomes
\begin{align*}
  \sum_{n \le N}  r_{\ell,2}(n) \frac{(N - n)^k}{\Gamma(k + 1)}
  &=
  \frac 1{2 \pi i}
  \int_{(a)}
    \mathcal{S}_{\ell}(z)
    \omega_2(z) e^{N z} z^{- k - 1} \, \dx z \\
  &\qquad+
  \Odi{\int_{(a)}
    \vert E(a,z,\ell)\vert  \, \vert e^{N z}\vert  \, \vert z\vert ^{- k - 1} \vert \omega_2(z)\vert  \, \vert \dx z\vert}.
\end{align*}
Using \eqref{omega-estim}-\eqref{z^-1}
and \eqref{expl-form-err-term-strong}, we see that the error term is
\begin{align*}
&\ll_{\ell}
 a^{- 1 / 2} e^{Na} 
 \Bigl( 
 \int_{-a}^{a} a^{{-k-1/2}} \dx y 
 +
 \int_{a}^{+\infty} y^{{-k-1/2}} \log^2(y/a)\,  \dx y 
 \Bigr)
 \\
&\ll_{k,\ell}
  e^{N a} a^{-k}
  \Bigl( 1+
 \int_{1}^{+\infty}v^{{-k-1/2}} \log^2 v \, \dx v
   \Bigr)
\ll_{k,\ell}
  e^{N a} a^{-k},
\end{align*}
provided that $k>1/2$.
Choosing $a = 1 / N$, the previous estimate becomes
\(
\ll_{k,\ell}
  N^{k}.
\)
Summing up, for  $k>  1/2$, we can write
\begin{equation}
\label{first-step-hl}
  \sum_{n \le N}  r_{\ell,2}(n) \frac{(N - n)^k}{\Gamma(k + 1)}
  =
  \frac 1{2 \pi i}
  \int_{(\frac{1}{N})}
    \mathcal{S}_{\ell}(z)
    \omega_2(z) e^{N z} z^{- k - 1} \, \dx z
  +
  \Odipg{k,\ell}{N^{k}}.
\end{equation}
We now insert \eqref{func-eq-omega} into \eqref{first-step-hl}, so that
the integral on the right-hand side of \eqref{first-step-hl} becomes
\begin{align}
\notag
  &
  \frac 1{2 \pi i}
  \int_{(\frac{1}{N})}
    \mathcal{S}_{\ell}(z)
    \Bigl( \frac 12 \Bigl( \frac \pi z \Bigr)^{1/2} - \frac12 \Bigr)
    e^{N z} z^{- k - 1} \, \dx z
  +
  \frac 1{2 \pi i}
  \int_{(\frac{1}{N})}
    \Bigl( \frac \pi z \Bigr)^{1/2}
    \mathcal{S}_{\ell}(z)
    \omega_2 \Bigl( \frac{\pi^2}z \Bigr)
     e^{N z} z^{- k - 1} \, \dx z \\
\label{hl-splitting}
  &=
  \I_1 + \I_2,
\end{align}
say. We now proceed to evaluate $\I_1$ and $\I_2$.
In the next two sections we will use \eqref{def-Sl} and obtain that
$I_1$ and $I_2$ split into a number of summands; in later sections we
will prove that we can exchange all summations and integrations, in
suitable ranges for $k$, using some properties of the non-trivial
zeros of the Riemann $\zeta$-function: see \S\ref{lemmas}.
Finally, we perform a change of variables that yields all the summands
in the statement of Theorem~\ref{Cesaro-gen-HL-average}.

\section{Evaluation of \texorpdfstring{$\I_1$}{I1}}
We need $k>  1/2$ in this section.
By a  direct computation we can write that
\begin{align*}
  \I_1
  &=
  \frac 1{4 \pi i}
   \frac{\Gamma(1/\ell)}{\ell}
  \int_{(\frac{1}{N})}
    \Bigl(  \frac{\pi^{1/2}}{z^{1/2}} - 1 \Bigr)
   e^{N z} z^{- k - 1-1/\ell}  \, \dx z
  -
  \frac {\pi^{1/2}}{4 \ell \pi i} 
  \int_{(\frac{1}{N})}
    \sum_{\rho} \Gamma\Bigl(\frac{\rho}{\ell}\Bigr) e^{N z} z^{- k -\rho/\ell - 3/2} \, \dx z \\
  &\qquad+
  \frac 1{4 \ell\pi i} 
  \int_{(\frac{1}{N})}
    \sum_{\rho}  \Gamma\Bigl(\frac{\rho}{\ell}\Bigr) e^{N z}  z^{- k - \rho/\ell - 1}  \, \dx z 
- 
  \frac {\log(2\pi)}{4\pi i} 
  \int_{(\frac{1}{N})}
 \Bigl(  \frac{\pi^{1/2}}{z^{1/2}} - 1 \Bigr)
  e^{N z}  z^{- k  - 1}  \, \dx z 
 \\& =
  \J_1 + \J_2 + \J_3+ \J_4,
\end{align*}
say. We see now how to evaluate $\J_1$, $\J_2$, $\J_3$ and $\J_4$.
The delicate point is the justification of the exchanges required to
deal with $\J_2$ and $\J_3$ (see \S\ref{exchange-rho-integral} for the
details), whereas the computations needed for $\J_1$ and $\J_4$ are
straightforward 
and immediately follow by using  the substitution $s=Nz$, 
by \eqref{Laplace-transf}. This way we get 
\begin{equation}
  \J_1  
 \label{J1-eval}
  =
  \frac{\pi^{1/2}}2    \frac{\Gamma(1/\ell)}{\ell}
 \frac{N^{k + 1/2+1/\ell}}{\Gamma(k + 3/2+1/\ell)}
  -
  \frac{\Gamma(1/\ell)}{2\ell}
   \frac{N^{k + 1/\ell}}{\Gamma(k + 1 + 1/\ell)}
\end{equation}
and
\begin{equation}
\J_4  
\label{J4-eval}
  =
- 
  \frac {\pi^{1/2}\log(2\pi)}{2\Gamma(k+3/2)} N^{k+1/2}
 +
  \frac {\log(2\pi)}{2 \Gamma(k+1)} N^{k}
  .
\end{equation}

\subsection{Evaluation of \texorpdfstring{$\J_2$}{J2}}

Exchanging the sum over $\rho$ with the integral (this can be done for $k>0$;
see \S\ref{exchange-rho-integral}) and using the substitution $s=Nz$, we have
\begin{align}
\notag
  \J_2
  &=
  -
  \frac{\pi^{1/2}}{2\ell}
  \sum_{\rho} \Gamma\Bigl(\frac{\rho}{\ell}\Bigr)
    \frac 1{2 \pi i}
    \int_{(\frac{1}{N})} e^{N z}  z^{-k -\rho/\ell - 3/2} \, \dx z \\
    \notag
  &=
  -
  \frac{\pi^{1/2}}{2\ell}
  \sum_{\rho} \Gamma\Bigl(\frac{\rho}{\ell}\Bigr)    
   N^{k + \rho/\ell + 1/2}
    \frac 1{2 \pi i}
    \int_{(1)} e^s s^{-k -\rho/\ell - 3/2}  \, \dx s   \\
    \label{J2-eval}
  &=
  -
  \frac{\pi^{1/2}}{2\ell}
  \sum_{\rho}   
  \frac{\Gamma(\rho/\ell)}
   {\Gamma(k + 3 / 2 + \rho/\ell)}
   N^{k + \rho/\ell + 1/2},
\end{align}
again by \eqref{Laplace-transf}.
By the Stirling formula \eqref{Stirling}, 
we remark that the series in $\J_2$ converges absolutely 
 for  $k>-1/2$. 

\subsection{Evaluation of \texorpdfstring{$\J_3$}{J3}}
Arguing as in \S\ref{exchange-rho-integral} with $-k-1$ which plays 
the role of $-k-3/2$ there, we see that we can exchange the sum with the integral 
provided that $k>1/2$. Hence, performing again the usual substitution $s=Nz$,
we can write
\begin{equation}
\label{J3-eval}
  \J_3 
  =
  \frac{1}{2\ell}
    \sum_{\rho}  \Gamma\Bigl(\frac{\rho}{\ell}\Bigr) N^{k + \rho/\ell}
    \frac 1{2 \pi i}
    \int_{(1)} e^s s^{- k - 1 - \rho/\ell} \, \dx s \\
  =
  \frac{1}{2\ell}
    \sum_{\rho} 
     \frac{\Gamma(\rho/\ell)}{\Gamma(k + 1 + \rho/\ell)} N^{k + \rho/\ell}.
\end{equation}
By the Stirling formula \eqref{Stirling}, 
we remark that the series in $\J_3$ converges absolutely 
 for  $k>0$. 

\section{Evaluation of \texorpdfstring{$\I_2$}{I2} and conclusion of the proof of Theorem \ref{Cesaro-gen-HL-average}}
 
 We need $k>  1$ in this section.
Using \eqref{def-Sl}  and the definition of $\omega_2 (\pi^2/ z )$ 
(see \eqref{Stilde-omega-def}) we have
\begin{align}
\notag
  \I_2
  &=
  \frac 1{2 \pi i} \frac{\Gamma(1/\ell)}{\ell}
  \int_{(\frac{1}{N})}
    \Bigl( \frac \pi z \Bigr)^{1/2}
    \Bigl( \sum_{j \ge 1} e^{- j^2 \pi^2 / z} \Bigr)
     e^{N z} z^{- k - 1- 1/\ell} \, \dx z \\
  &
 \notag
  \qquad
  - \frac 1{2 \ell \pi i}
  \int_{(\frac{1}{N})}
    \Bigl( \frac \pi z \Bigr)^{1/2}
    \Bigl( \sum_{j \ge 1} e^{- j^2 \pi^2 / z} \Bigr)
    \Bigl( \sum_{\rho} z^{-\rho/\ell}  \Gamma\Bigl(\frac{\rho}{\ell}\Bigr) \Bigr)
     e^{N z} z^{- k - 1} \, \dx z 
\\
 \label{I2-splitting}
&
-    \frac {\log(2\pi)}{2 \pi i}
  \int_{(\frac{1}{N})}
    \Bigl( \frac \pi z \Bigr)^{1/2}
    \Bigl( \sum_{j \ge 1} e^{- j^2 \pi^2 / z} \Bigr) 
     e^{N z} z^{- k - 1} \, \dx z 
  =   \J_5   +   \J_6    +  \J_7,
  \end{align}
say. We see now how to evaluate  $ \J_5, \J_6$ and $\J_7$.
The proof in this section is more delicate than in the previous one.
We have to justify inversion as before, but we are then faced with the
problem of dealing with series containing values of the Bessel
functions, arising from the ``modular'' terms.
We refer to \S\ref{sums-abs-conv} for a detailed discussion of the
problem.

\subsection{Evaluation of \texorpdfstring{$\J_5$}{J5}}

By means of the substitution $s = N z$, since
the exchange  is justified   in \S\ref{exchange-ell-integral} for $k> 1/2 - 1/\ell$, we get
\begin{equation*}
  \J_5 
  =
  \pi^{1 / 2} \frac{\Gamma(1/\ell)}{\ell}
  N^{k + 1 / 2 + 1/\ell}
  \sum_{j \ge 1}  
    \frac 1{2 \pi i}
    \int_{(1)}
      e^s e^{- j^2 \pi^2 N / s} s^{- k - 3 / 2 - 1 /\ell} \, \dx s.
\end{equation*}
Setting $u = 2 \pi j N^{1/2}$ in \eqref{Bessel-def}, we obtain
\begin{equation}
\label{J-nu}
  J_\nu \bigl( 2 \pi j N^{1/2} \bigr)
  =
  \frac{(\pi j N^{1/2})^\nu}{2 \pi i}
  \int_{(1)}  e^s e^{- j^2 \pi^2 N / s} s^{- \nu -1}\, \dx s,
\end{equation}
and  hence we have
\begin{equation}
\label{J5-eval} 
 \J_5
=
  \frac{N^{k / 2 + 1 / 4 + 1/(2\ell)}}{\pi^{k + 1/\ell}}
  \frac{\Gamma(1/\ell)}{\ell}
  \sum_{j \ge 1}  
    \frac{J_{k + 1 / 2 + 1/ \ell} (2 \pi j N^{1 / 2})}{j^{k + 1 / 2 + 1/ \ell}}.
\end{equation}
The absolute convergence of the series in $ \J_5$  is studied
in \S\ref{sums-abs-conv}.

\subsection{Evaluation of \texorpdfstring{$\J_6$}{J6}}

With the same substitution used before, since
the double exchange between sums and the line integral
is justified in \S\ref{exchange-double-sum-ell-rho}
for $k>1$, we see that
\begin{equation*}
  \J_6
  =
  - \frac{\pi^{1 / 2}}{\ell}
  \sum_{\rho}
         \Gamma\Bigl(\frac{\rho}{\ell}\Bigr)  
        N^{k + 1 / 2 + \rho/\ell}
        \sum_{j \ge 1}  
      \Bigl(
        \frac 1{2 \pi i}
        \int_{(1)}
          e^s e^{- j^2 \pi^2 N / s} s^{- k - 3 / 2 - \rho/\ell} \, \dx s
      \Bigr).
\end{equation*}
Using  \eqref{J-nu}, we get
\begin{equation}
  \J_6
  =
\label{J6-eval}
  -
  \frac{N^{k / 2 + 1 / 4}}{\pi^k}
  \sum_{\rho} \Gamma\Bigl(\frac{\rho}{\ell}\Bigr)  \frac{N^{\rho / (2\ell)}}{\pi^{\rho/\ell}}
    \sum_{j \ge 1}  
      \frac{J_{k + 1 / 2 + \rho/\ell} (2 \pi j N^{1 / 2})}
           {j^{k + 1 / 2 + \rho/\ell}}.
\end{equation}
In this case, the absolute convergence of the series in $\J_{6}$ 
is more delicate; such a treatment is again described 
in \S\ref{sums-abs-conv}.
 
\subsection{Evaluation of \texorpdfstring{$\J_7$}{J7}}

With the same substitution used before, since
the  exchange between sum and the line integral
is justified in \S\ref{exchange-ell-integral} for $k> 1/2$, we see that
\begin{equation*}
  \J_7
  =
  -\pi^{1/2}\log(2\pi) N^{k+1/2}
\sum_{j \ge 1} 
\frac {1}{2 \pi i}
  \int_{(1)} 
     e^{s} e^{- j^2 \pi^2 N / s}   s^{- k - 3/2} \, \dx s 
\end{equation*}
Using  \eqref{J-nu}, we get
\begin{equation}
  \J_7
  =
\label{J7-eval}
 -  \frac{\log(2\pi)}{\pi^{k}} N^{k/2+1/4}
     \sum_{j \ge 1}  
      \frac{J_{k + 1 / 2} (2 \pi j N^{1 / 2})}
           {j^{k + 1 / 2}}.
\end{equation}
The absolute convergence of the series in $\J_{7}$  
 is studied
in \S\ref{sums-abs-conv}.

Finally, inserting \eqref{J1-eval}--\eqref{J7-eval} into
\eqref{hl-splitting} and \eqref{first-step-hl}, we obtain
\begin{equation}
\label{expl-form-HL-bis}
  \sum_{n \le N} r_{\ell,2}(n) \frac{(N - n)^k}{\Gamma(k + 1)}
  =
  N^k \sum_{j=1}^7 \M_{j,\ell,k}(N)
  +
  \Odipg{k,\ell}{N^{k}},
\end{equation}
for $k > 1$. 
Theorem \ref{Cesaro-gen-HL-average} follows dividing
\eqref{expl-form-HL-bis} by $N^{k}$.

\section{Lemmas}
\label{lemmas}

We recall some basic facts in complex analysis.
First, if $z = a + i y$ with $a > 0$, we see that for complex $w$ we
have
\begin{align*}
  z^{-w}
  &=
  \vert z \vert^{-w} \exp( - i w \arctan(y / a)) \\
  &=
  \vert z \vert^{-\Re(w) - i \Im(w)} \exp( (- i \Re(w) + \Im(w)) \arctan(y / a)),
\end{align*}
so that
\begin{equation}
\label{z^w}
  \vert z^{-w} \vert
  =
  \vert z \vert^{-\Re(w)} \exp(\Im(w) \arctan(y / a)).
\end{equation}
We also recall that, uniformly for $x \in [x_1, x_2]$, with $x_1$ and
$x_2$ fixed, and for $|y| \to +\infty$, by the Stirling formula we have
\begin{equation}
\label{Stirling}
  \vert \Gamma(x + i y) \vert
  \sim
  \sqrt{2 \pi}
  e^{- \pi |y| / 2} |y|^{x - 1 / 2};
\end{equation}
see, e.g., Titchmarsh \cite[\S4.42]{Titchmarsh1988}.

We will need the  following lemmas from Languasco and Zaccagnini \cite{LanguascoZ2017c}.

\begin{Lemma}[See Lemma 1 of \cite{LanguascoZ2017c}] 
\label{Linnik-lemma-gen}   
Let $\ell\ge 1$ be an integer,  $z = a + iy$, where $a > 0$ and $y \in \R$
and let $\mathcal{S}_{\ell}(z)$ be defined as in \eqref{def-Sl}.
Then $\Stilde_{\ell}(z) = \mathcal{S}_{\ell}(z) + E(a,z,\ell)$ where
\begin{equation}
\label{expl-form-err-term-strong}
  E(a,z,\ell)
  \ll_\ell
  \vert z \vert^{1/2}
  \begin{cases}
    1 & \text{if $\vert y \vert \leq a$} \\
    1 +\log^2 (\vert y\vert/a) & \text{if $\vert y \vert > a$.}
  \end{cases}
\end{equation}
\end{Lemma}

\begin{Lemma}[See Lemma 2 of \cite{LanguascoZ2017c}] 
\label{series-int-zeros} 
Let $\ell\ge 1$ be an integer, let $\rho = \beta + i \gamma$ run over the non-trivial zeros of the Riemann
zeta-function and $\alpha > 1$ be a parameter.
The series
\[
  \sum_{\rho \colon \gamma > 0}
  \Bigl(\frac{\gamma}{\ell}\Bigr)^{\beta/\ell-1/2}
    \int_1^{+\infty} \exp\Bigl( - \frac{\gamma}{\ell} \arctan\frac 1u \Bigr)
      \frac{\dx u}{u^{\alpha+\beta/\ell}}
\]
converges provided that $\alpha > 3/2$.
For $\alpha \le 3/2$ the series does not converge.
The result remains true if we insert in the integral a factor
$(\log u)^c$, for any fixed $c \ge 0$.
\end{Lemma}
\begin{Lemma}[See Lemma 3 of \cite{LanguascoZ2017c}]  
\label{series-int-zeros-alt-sign}
Let $\ell\ge 1$ be an integer,   $\alpha > 1$, $z=a+iy$, $a\in(0,1)$ and $y\in \R$.
Let further $\rho = \beta+i\gamma$ run over the non-trivial zeros of
the Riemann zeta-function.  We have
\[
  \sum_{\rho}
    \Bigl\vert \frac{\gamma}{\ell}\Bigr\vert ^{\beta/\ell-1/2}
    \int_{\Y_1 \cup \Y_2} \exp\Bigl(\frac{\gamma }{\ell} \arctan\frac{y}{a} - \frac\pi2\Bigl \vert \frac{\gamma}{\ell}  \Bigr\vert\Bigr)
      \frac{\dx y}{\vert z \vert ^{\alpha+\beta/\ell}}
  \ll_{\alpha,\ell}
  a^{1-\alpha-1/\ell},
\]
where $\Y_1=\{y\in \R\colon y\gamma \leq 0\}$ and 
$\Y_2=\{y\in [-a,a] \colon y\gamma > 0\}$.
The result remains true if we insert in the integral a factor
$(\log (\vert y\vert /a))^c$, for any fixed $c \ge 0$.
\end{Lemma}

\section{Interchange of the series over zeros with the line integral in \texorpdfstring{$\J_2, \J_3$}{J2, J3}} 
\label{exchange-rho-integral}

We need $k>1/2$ in this section. For $\J_2$ we have to establish the convergence of
\begin{equation}
\label{conv-integral-J2-J3}
  \sum_{\rho}
    \Bigl\vert \Gamma\Bigl(\frac{\rho}{\ell}\Bigr) \Bigr\vert 
     \int_{(\frac{1}{N})} \vert e^{N z}\vert  \vert z\vert ^{- k - 3/2} \vert z^{- \rho/\ell}\vert  \, \vert \dx z\vert,
\end{equation}
where, as usual, $\rho=\beta + i \gamma$ runs over the non-trivial 
zeros of the Riemann zeta-function.
By \eqref{z^w} and the Stirling formula \eqref{Stirling}, we are left
with estimating
\begin{equation}
\label{conv-integral-J2-J3-1}
  \sum_{\rho}
    \Bigl\vert \frac{\gamma}{\ell}\Bigr\vert ^{\beta/\ell-1/2}
    \int_{\R}  \exp\Bigl(\frac{\gamma }{\ell} \arctan(N y) -\frac\pi2\Bigl \vert \frac{\gamma}{\ell} \Bigr\vert \Bigr)
             \frac{\dx y}{\vert z \vert ^{k + 3/2 +\beta/\ell}}. 
\end{equation}
We have just to consider the case $\gamma y >0$, $\vert y \vert > 1/N$
since in the other cases the total contribution is  $\ll_k N^{k + 1/2+1/\ell }$
by Lemma \ref{series-int-zeros-alt-sign} with $\alpha=k+3/2$ and $a=1/N$.
By symmetry, we may assume that $\gamma > 0$.
We have that the integral in \eqref{conv-integral-J2-J3-1} is
\begin{align*}
    &
  \ll_{\ell}
    \sum_{\rho}
    \Bigl\vert \frac{\gamma}{\ell}\Bigr\vert ^{\beta/\ell-1/2}
    \int_{1 / N}^{+\infty} \exp\Bigl( - \frac{\gamma }{\ell} \arctan\frac 1{N y} \Bigr)
      \frac{\dx y}{y^{k + 3/2 +\beta/\ell}} \\
  &=
  N^{k+1/2}
  \sum_{\rho \colon \gamma > 0}
    N^{\beta/\ell}
     \Bigl( \frac{\gamma}{\ell}\Bigr) ^{\beta/\ell-1/2}
    \int_1^{+\infty} \exp\Bigl( - \frac{\gamma }{\ell} \arctan\frac 1u \Bigr)
      \frac{\dx u}{u^{k + 3/2 +\beta/\ell}}.
\end{align*}
For $k > 0$, this is $\ll_{k,\ell} N^{k + 1/2 +1/\ell}$ by
Lemma~\ref{series-int-zeros}. This implies
that the integrals in \eqref{conv-integral-J2-J3-1} and in \eqref{conv-integral-J2-J3} 
are both  $\ll_{k,\ell} N^{k +1/2+ 1/\ell}$, and
hence this exchange step for $\J_2$ is fully justified.

For $\J_3$, we have to consider
\begin{equation}
\label{conv-integral-J3}
  \sum_{\rho}
    \Bigl\vert \Gamma\Bigl(\frac{\rho}{\ell}\Bigr) \Bigr\vert 
     \int_{(\frac{1}{N})} \vert e^{N z}\vert  \vert z\vert ^{- k - 1} \vert z^{- \rho/\ell}\vert  \, \vert \dx z\vert.
\end{equation}
We can repeat the same reasoning we used for $\J_2$ just replacing $k+3/2$ 
with $k+1$. This means that we need $k>1/2$ here 
to get that the  integral in \eqref{conv-integral-J3} is
$\ll_{k,\ell} N^{k + 1/\ell}$, and that
  this exchange step for $\J_3$ is fully justified too.

\section{Interchange of the series over \texorpdfstring{$j$}{j} with the line integral \texorpdfstring{in $\J_5, \J_7$}{}} 
\label{exchange-ell-integral}

We need $k>  1/2  $ in this section. For $\J_5$
we  have to  establish the convergence of
\begin{equation}
\label{conv-integral-J5}
  \sum_{j \ge 1}  
    \int_{(\frac{1}{N})} 
    \vert e^{N z}\vert  \vert z\vert ^{- k - 3/2 - 1/\ell}
    e^{-\pi^2 j^2 \Re(1/z)}  \, \vert \dx z\vert .
\end{equation}
A trivial computation gives 
\begin{equation}
\label{real-part-estim}
\Re(1/z)
=
\frac{N}{1+N^2y^2}
\gg
\begin{cases}
N & \text{if}\ \vert y \vert \leq 1/N, \\
1/(Ny^2) & \text{if}\ \vert y \vert > 1/N.
\end{cases}
\end{equation}
By \eqref{real-part-estim}, we can write that
the  quantity in \eqref{conv-integral-J5} is
\begin{equation} 
    \label{J4-split}
\ll_\ell 
  \sum_{j \ge 1}  
    \int_{0}^{1/N}  
    \frac{e^{- j^2 N}}{\vert z\vert ^{k + 3/2 + 1/\ell}}  \, \dx y
    +
 \sum_{j \ge 1}  
    \int_{1/N}^{+\infty} 
    \frac{e^{- j^2 / (Ny^2)}}{\vert z\vert ^{k + 3/2 + 1/\ell}}   \, \dx y
    =
U_1+U_2,
\end{equation}
say, since the $\pi^2$ factor in the exponential function
is negligible. Using \eqref{omega-estim}-\eqref{z^-1},  we have 
\begin{equation}
\label{U1-estim}
U_1
\ll_\ell
 N^{k+ 1/2 + 1/\ell}
\omega_{2} (N) 
 \ll_\ell  
 N^{k+1/\ell}
\end{equation}
and
\begin{align}
\notag
U_2&
\ll_\ell
\sum_{j \ge 1}
    \int_{1/N}^{+\infty} 
    \frac{e^{- j^2 / (Ny^2)}}{y^{k + 3/2 + 1/\ell}}  \, \dx y
  \ll_\ell
N^{k/2+1/4+1/(2\ell)} \sum_{j \ge 1}  
\frac{1}{j^{k+ 1/2 + 1/\ell}}
    \int_{0}^{j^2 N} u^{ k/2 - 3/4 +1/(2\ell)} 
    e^{- u}  \, \dx u 
\\
\label{U2-estim}
&\leq
\Gamma \Bigl( \frac{2k+1+2\ell}{4} \Bigr)
N^{k/2+1/4+1/(2\ell)}  \sum_{j \ge 1}  
\frac{1}{j^{k+1/2+1/\ell}}
\ll_{k,\ell}
N^{k/2+1/4+1/(2\ell)},
\end{align}
provided that $k>1/2-1/\ell$,
where  we used the substitution
$u=j^2 / (Ny^2)$.
Inserting \eqref{U1-estim}-\eqref{U2-estim}
into \eqref{J4-split} we get, for $k>1/2-1/\ell$, that
the quantity in \eqref{conv-integral-J5}
is $\ll N^{k+1/\ell}$ and so it is for  $\J_5$.

For $\J_7$
we  have to  establish the convergence of
\begin{equation}
\label{conv-integral-J7}
  \sum_{j \ge 1}  
    \int_{(\frac{1}{N})} 
    \vert e^{N z}\vert  \vert z\vert ^{- k - 3/2}
    e^{-\pi^2 j^2 \Re(1/z)}  \, \vert \dx z\vert .
\end{equation}
We can repeat the same reasoning we used for $\J_5$ just replacing $k+3/2+1/\ell$ 
with $k+3/2$. This means that we need $k>1/2$ here 
to get that the  integral in \eqref{conv-integral-J7} is
$\ll_{k,\ell} N^{k }$, and that
  this exchange step for $\J_7$ is fully justified too.

\section{Interchange of series with the line integral in \texorpdfstring{$\J_6$}{J6}}
\label{exchange-double-sum-ell-rho}

We need $k>1$ in this section.
We first have to  establish the convergence of
\begin{equation}
\label{conv-integral-J6}
  \sum_{j \ge 1}  
    \int_{(\frac{1}{N})} 
         \Bigl\vert \sum_{\rho} \Gamma\Bigl(\frac{\rho}{\ell}\Bigr)z^{- \rho/\ell} \Bigr\vert
    \vert e^{N z}\vert  \vert z\vert ^{- k - 3/2}
    e^{-\pi^2 j^2 \Re(1/z)}  \, \vert \dx z\vert .
\end{equation}
Using the Prime Number Theorem and 
\eqref{expl-form-err-term-strong}, 
we first remark that
\begin{equation} 
\label{sum-over-rho-new-HL} 
 \Bigl\vert \sum_{\rho} \Gamma\Bigl(\frac{\rho}{\ell}\Bigr)z^{- \rho/\ell} \Bigr\vert
  \ll_\ell
  N^{1/\ell} + \vert z\vert ^{1/2}
    \log^2 (2N\vert y\vert) . 
\end{equation}

By \eqref{real-part-estim} 
and \eqref{sum-over-rho-new-HL}, we can write that
the  quantity in \eqref{conv-integral-J6} is
\begin{align}
\notag 
\ll_\ell N^{1/\ell}
  \sum_{j \ge 1}  
    \int_{0}^{1/N}  
    \frac{e^{- j^2 N}}{\vert z\vert ^{k + 3/2}}  \,&  \dx y
    +
 N^{1/\ell}\sum_{j \ge 1}
    \int_{1/N}^{+\infty}   
    \frac{e^{- j^2 / (Ny^2)} }{\vert z\vert ^{k + 3/2}}  \, \dx y
    \\ 
    \label{HL-split}
    &
    +
    \sum_{j \ge 1}  
    \int_{1/N}^{+\infty}   \log^{2}(2Ny)
    \frac{e^{- j^2 / (Ny^2)}}{\vert z\vert ^{k + 1}}   \, \dx y
    =
V_1+V_2+V_{3},
\end{align}
say.
$V_{1}$  can be estimated exactly as $U_{1}$
in Section \ref{exchange-ell-integral} 
and we get $V_{1} \ll_{k,\ell} N^{k+ 1/\ell}$.
For $V_{2}$ we can work analogously to $U_{2}$ thus obtaining
\begin{align}
\notag
V_2&
\ll_{k,\ell}
  N^{1/\ell} \sum_{j \ge 1}  
    \int_{1/N}^{+\infty} 
    \frac{e^{- j^2 / (Ny^2)}}{y^{k + 3/2}}  \, \dx y
\ll_{k,\ell}
  N^{k/2+1/4+1/\ell} \sum_{j \ge 1}  
\frac{1}{j^{k+1/2}}
    \int_{0}^{j^2 N} u^{ k/2 - 3/4} 
    e^{- u}  \, \dx u 
\\ 
\notag
&\ll_{k,\ell}
\Gamma \Bigl( \frac{2k+1}{4} \Bigr)
N^{k/2+1/4+1/\ell}  \sum_{j \ge 1}  
\frac{1}{j^{k+1/2}}
\ll_{k,\ell}
N^{k/2+1/4+1/\ell},
\end{align}
provided that $k>1/2$,
where  we used the substitution
$u=j^2 / (Ny^2)$.
Hence, we have 
\begin{equation}
\label{V1-V2-estim}
V_{1} + V_{2} \ll_{k,\ell} N^{k+1/\ell},
\end{equation}
provided that  $k>  1/2$.

Using the substitution $u=j^2 / (Ny^2)$,
we obtain
\begin{align}
\notag
V_3
&\ll_{k,\ell}
 \sum_{j \ge 1}  
    \int_{1/N}^{+\infty} \log^{2}(2Ny) 
    \frac{e^{- j^2 / (Ny^2)}}{y^{k + 1}}   \, \dx y
=
\frac{N^{k/2}}{8} \sum_{j \ge 1}  
\frac{1}{j^k}
    \int_{0}^{j^2 N} u^{ k/2 - 1} 
    \log^{2}\Bigl( \frac{4j^{2}N}{u} \Bigr)
    e^{- u}  \, \dx u.
\end{align}
Hence, a direct computation shows that
\begin{align}
\notag
V_3
& \ll_{k,\ell}
N^{k/2} \sum_{j \ge 1}  
\frac{\log^{2}(j N)}{j^k}
    \int_{0}^{j^2 N} u^{ k/2 - 1}     e^{- u}  \, \dx u
    +
N^{k/2} \sum_{j \ge 1}  
\frac{1}{j^k}
     \int_{0}^{j^2 N} u^{ k/2 - 1}  \log^{2} (u)\,
    e^{- u}  \, \dx u
\\
\label{V3-estim}
&\ll_{k,\ell}
\Gamma(k/2)
N^{k/2} \sum_{j \ge 1}  
\frac{\log^{2}(j N)}{j^k}
    +
N^{k/2}
\ll_{k,\ell}
N^{k/2}\log^{2} N
\end{align}
provided that $k>1$.
Inserting \eqref{V1-V2-estim}-\eqref{V3-estim}
into \eqref{HL-split} we get, for $k>1$, that
the quantity in \eqref{conv-integral-J6}
is $\ll_{k,\ell} N^{k+1/\ell}$.

Now we have to  establish the convergence of
\begin{equation}
\label{conv-integral-5}
  \sum_{j \ge 1}  
     \sum_{\rho} \Bigl\vert \Gamma\Bigl(\frac{\rho}{\ell}\Bigr) \Bigr\vert 
    \int_{(\frac{1}{N})} 
   \vert e^{N z}\vert  \vert z\vert ^{- k - 3/2}
         \vert z^{- \rho/\ell} \vert
    e^{-\pi^2 j^2 \Re(1/z)}  \, \vert \dx z\vert .
\end{equation}
By symmetry, we may assume that $\gamma > 0$.
For $y \in(-\infty, 0]$ we have
$\gamma \arctan(y/a) -\frac \pi2  \gamma  \le - \frac \pi2  \gamma$.
Using \eqref{real-part-estim}, \eqref{z^-1} and the Stirling formula \eqref{Stirling}, 
the quantity we are estimating becomes
\begin{align}
\notag
&\ll
  \sum_{j \ge 1}  
  \sum_{\rho \colon \gamma>0}
 \Bigl( \frac{\gamma}{\ell}\Bigr)^{\beta/\ell-1/2}
   \exp\Bigl(  -\frac \pi2  \frac{\gamma}{\ell}\Bigr)
   \Bigl(
              \int_{-1/N}^{0}  
         N^{k + 3/2 + \beta/\ell}\ e^{- j^2 N}  \, \dx y
        +
        \int_{-\infty}^{-1/N}
         \frac{e^{- j^2 / (Ny^2)}}{\vert y \vert^{k + 3/2 + \beta/\ell}}  \, \dx y
        \Bigr)
    \\
    \notag
   & \ll_{k,\ell}
N^{k+1/2+1/\ell}
 \sum_{j \ge 1}      e^{- j^2 N} 
\sum_{\rho \colon \gamma>0}
 \Bigl(\frac{\gamma}{\ell}\Bigr)^{\beta/\ell-1/2}
\exp\Bigl(  -\frac \pi2  \frac{\gamma}{\ell}\Bigr)
\\
\notag
&\hskip1cm
+
N^{k/2+1/4} \sum_{j \ge 1}  
\frac{1}{j^{k+1/2}}
\sum_{\rho \colon \gamma>0}
\frac{N^{\beta/(2\ell)}}{j^{\beta/\ell}}
 \Bigl(\frac{\gamma}{\ell}\Bigr)^{\beta/\ell-1/2}
\exp\Bigl(  -\frac \pi2 \frac{\gamma}{\ell}\Bigr)
    \int_{0}^{j^2 N} u^{ k/2 - 3/4 + \beta/(2\ell)} 
    e^{- u}  \, \dx u 
    \\
    \notag
    &
    \ll_{k,\ell}
     N^{k+1/\ell}
     + 
\Bigl(
  \max_{0\le b \le 1} \Gamma\Bigl( \frac{b}{2\ell} +\frac{k}2+ \frac14\Bigr) 
\Bigr) 
N^{\frac{k}{2}+\frac14+\frac{1}{2\ell}} \sum_{j \ge 1}  
\frac{1}{j^{k+1/2}}
\sum_{\rho \colon \gamma>0} 
 \Bigl(\frac{\gamma}{\ell}\Bigr)^{\beta/\ell-1/2}
\exp\Bigl(  -\frac \pi2 \frac{\gamma}{\ell}\Bigr)
    \\
    \label{easy-case}
    & \ll_{k,\ell}  
 N^{k+1/\ell},
\end{align}
provided that $k>1/2$, 
where  we used the substitution $u = - j^2 / (Ny^2)$, 
 \eqref{omega-estim} and standard density estimates.

Let now $y>0$.
Using the Stirling formula \eqref{Stirling} and \eqref{real-part-estim}, we can write that
the  quantity in \eqref{conv-integral-5} is
\begin{align}
\notag
\ll_{k,\ell} &
  \sum_{j \ge 1}  
  \sum_{\rho \colon \gamma>0}
    \Bigl(\frac{\gamma}{\ell}\Bigr)^{\beta/\ell-1/2}
\exp\Bigl(  -\frac \pi4 \frac{\gamma}{\ell}\Bigr)
              \int_{0}^{1/N}  
        \frac{e^{- j^2 N}}{\vert z\vert ^{k + 3/2 + \beta/\ell}}  \, \dx y
   \\ 
    \label{HL-split-1}
   &+
\sum_{j \ge 1}
  \sum_{\rho \colon \gamma>0}
   \Bigl(\frac{\gamma}{\ell}\Bigr)^{\beta/\ell-1/2}
      \int_{1/N}^{+\infty}  
\exp\Bigl( \frac{\gamma}{\ell} (\arctan(N y) -\frac \pi2) \Bigr)  
    \frac{e^{- j^2 / (Ny^2)}}{\vert z\vert ^{k + 3/2 + \beta/\ell}}  \, \dx y
    =
W_1+W_2,
\end{align}
say.
Using \eqref{z^-1} and \eqref{omega-estim},  we have that
\begin{equation}
\label{W1-estim}
W_1
\ll_{k,\ell}
N^{k+1/2+1/\ell}
 \sum_{j \ge 1}      e^{- j^2 N} 
\sum_{\rho \colon \gamma>0}
\Bigl(\frac{\gamma}{\ell}\Bigr)^{\beta/\ell-1/2}
\exp\Bigl(  -\frac \pi4 \frac{\gamma}{\ell}\Bigr)
\ll_{k,\ell}
N^{k +1/\ell},
\end{equation}
by standard density estimates.
Moreover, we get
\begin{align*}
W_2
&\ll_{k,\ell}
\sum_{j \ge 1}
\sum_{\rho \colon \gamma>0}
\Bigl(\frac{\gamma}{\ell}\Bigr)^{\beta/\ell-1/2}
        \int_{1/N}^{+\infty} y^{- k - 3/2 - \beta/\ell}  
\exp\Bigl( - \frac{\gamma}{\ell N y} - \frac{j^2}{Ny^2} \Bigr)  \, \dx y
\\
&\ll_{k,\ell}
N^{k/2+1/4} 
\sum_{j \ge 1}
\frac{1}{j^{k+1/2}}
\sum_{\rho \colon \gamma>0}
\frac{N^{\beta/(2\ell)} \gamma^{\beta/\ell-1/2}} {j^{\beta/\ell}}
    \int_{0}^{j \sqrt{N}} v^{ k - 1/2 + \beta/\ell} 
    \exp\Bigl(-\frac{\gamma v}{\ell j\sqrt{N}}- v^2\Bigr) 
    \, \dx v,
\end{align*}
in which we used the substitution
$v^2=j^2 / (Ny^2)$.
We remark that, for $k>1$, 
we can set $\eps=\eps(k)=(k-1)/2>0$ and that 
$k -\eps =(k+1)/2>1$.
We further remark that $\max_{v} (v^{k -\eps}e^{-v^2})$
is attained at $v_0=((k -\eps)/2)^{1/2}$, and hence we obtain,
for $N$ sufficiently large, that
\[
W_2
\ll_{k,\ell}
N^{k/2+1/4}
  \sum_{j \ge 1}  
\frac{1}{j^{k+1/2}}
\sum_{\rho \colon \gamma>0}
\frac{N^{\beta/(2\ell)} \gamma^{\beta/\ell-1/2}} {j^{\beta/\ell}}
    \int_{0}^{j \sqrt{N}}  
    v^{\beta/\ell-1/2+\eps}\exp\Bigl(-\frac{\gamma v}{\ell j\sqrt{N}}\Bigr)
     \, \dx v .
\]
Making the substitution $u= \gamma v/(j \sqrt{N})$,  we have
\begin{align}
\notag
W_2
&\ll_{k,\ell}
N^{k/2+1/2+\eps/2}
  \sum_{j \ge 1}  
\frac{1}{j^{k -\eps}}
\sum_{\rho \colon \gamma>0}
\frac{N^{\beta/\ell}} {\gamma^{1+\eps}}
    \int_{0}^{\gamma}  
    u^{\beta/\ell-1/2+\eps} e^{- u}
     \, \dx u\\
\label{W2-estim}
&\ll_{k,\ell}
N^{\frac{k}{2}+\frac12+\frac{1}{\ell}+\frac{\eps}{2}}
  \sum_{j \ge 1}  
\frac{1}{j^{k -\eps}}
\sum_{\rho \colon \gamma>0}
\frac{ 1 } {\gamma^{1+\eps}}
\Bigl(
\max_{0\le b \le 1} \Gamma \Bigl( \frac{b}{\ell} + \frac12+\eps \Bigr)  
\Bigr)  
\ll_{k,\ell}
 N^{\frac{3}{4}k+\frac14+\frac{1}{\ell}},
\end{align} 
by standard density estimates 
 and the definition of $\eps$.
Inserting \eqref{W1-estim}-\eqref{W2-estim}
into \eqref{HL-split-1} and recalling \eqref{easy-case},
we get, for $k>1$, that
the quantity in \eqref{conv-integral-5}
is $\ll_{k,\ell} N^{k+1/\ell}$.

\section{Absolute convergence of \texorpdfstring{$\J_{5}, \J_{6}$ and $\J_{7}$}{J5, J6 and J7}}
\label{sums-abs-conv}

\medskip 
Using, for $\nu>0$ fixed, $u\in \R$ and $u\to+\infty$, the estimate
\begin{equation}
\label{poisson-estim}
\vert J_\nu(u) \vert
\ll_\nu 
 u^{-1/2}
\end{equation}
which immediately follows from \eqref{Lebedev-asymp}, and performing
 a direct computation, we obtain that $ \J_5$
converges absolutely for $k> -1/\ell$ (and for $N$ sufficiently large)
and that $\J_5\ll_{k,\ell} N^{k/2+1/(2\ell)}$.
 
Again using \eqref{poisson-estim} and performing
 a direct computation as in the previous case, we obtain that $ \J_7$
converges absolutely for $k> 0$ (and for $N$ sufficiently large)
and that $\J_7\ll_{k,\ell} N^{k/2}$.
 
For the study of the absolute convergence of the series in $\J_6$ we have 
a different situation. 
In this case it is better to come back
to the Sonine 
representation of the Bessel functions \eqref{Bessel-def} on the line
$\Re(s)=1$. Using the usual substitution $s=Nz$,
we are led to consider the quantity
\begin{align*}
\sum_{\rho} 
 \Bigl\vert
  \Gamma\Bigl(\frac{\rho}{\ell}\Bigr) \frac{N^{\rho / (2\ell)}}{\pi^{\rho/\ell}} 
  \Bigr\vert
  &  \sum_{j \ge 1}  
       \Bigl\vert
       \frac{ J_{k + 1 / 2 + \rho/\ell} (2 \pi j N^{1 / 2}) }
           {j^{k + 1 / 2 + \rho/\ell}} 
           \Bigr\vert
           \\
&           
    \ll_{k}
    N^{-k/2-1/4}
     \sum_{\rho} \Bigl\vert \Gamma\Bigl(\frac{\rho}{\ell}\Bigr)  \Bigr\vert 
     \sum_{j \ge 1}    
    \int_{(\frac{1}{N})} 
   \vert e^{N z}\vert  \vert z\vert ^{- k - 3/2}
         \vert z^{- \rho/\ell} \vert
    e^{-\pi^2 j^2 \Re(1/z)}  \, \vert \dx z\vert,
\end{align*}
which is very similar to the one in \eqref{conv-integral-5}; the only difference is that
the sums are interchanged.
The argument used in \eqref{conv-integral-5}-\eqref{W2-estim} can be applied in this case
too thus showing that the double series in $\J_{6}$ converges absolutely for 
$k>1$.

\bigskip
\noindent
We thank the Referee for a very careful reading of the first version
of this paper.

\renewcommand{\bibliofont}{\normalsize}

\vskip1cm 
\noindent
\begin{tabular}{l@{\hskip 20mm}l}
Alessandro Languasco               & Alessandro Zaccagnini\\
Universit\`a di Padova     & Universit\`a di Parma\\
 Dipartimento di Matematica  & Dipartimento di Matematica, \\
 ``Tullio Levi-Civita'' &  Fisica e Informatica\\\
Via Trieste 63                & Parco Area delle Scienze, 53/a \\
35121 Padova, Italy            & 43124 Parma, Italy\\
{\it e-mail}: alessandro.languasco@unipd.it        & {\it e-mail}:
alessandro.zaccagnini@unipr.it  
\end{tabular}

\end{document}